\documentclass[12pt,leqno]{article}

\usepackage{amssymb,amsmath,amsthm}

\newtheorem{theorem}{Theorem}[section]

\newtheorem{proposition}[theorem]{Proposition}

\newtheorem{corollary}[theorem]{Corollary}

\theoremstyle{definition}
\newtheorem{definition}[theorem]{Definition}

\theoremstyle{remark}
\newtheorem{remark}[theorem]{Remark}

\numberwithin{equation}{section}

\newcommand{\wedgeco}%
{\displaystyle\operatornamewithlimits{\wedge}_{\raise3pt\hbox{,}}}

\newcommand{\ostar}{\odot\kern-6.4pt\ast}

\def\zs#1{_{\lower2pt\hbox{$\scriptstyle#1$}}}

\newcommand{\hb}{\hbar}

\newcommand{\pa}{\partial}

\newcommand{\wh}{\widehat}

\newcommand{\od}{\overset{\text{\rm def}}{=}}

\newcommand{\const}{\operatorname{const}}
\newcommand{\Exp}{\operatorname{Exp}}

\newcommand{\diag}{\operatorname{diag}}

\newcommand{\id}{\operatorname{id}}
\newcommand{\ad}{\operatorname{ad}}
\newcommand{\spp}{\operatorname{sp}}

\newcommand{\bS}{\mathbf{S}}

\newcommand{\bR}{\mathbb{R}}

\newcommand{\bH}{\mathbb{H}} 
\newcommand{\bC}{\mathbb{C}}

\newcommand{\cH}{\mathcal{H}} 
 
\newcommand{\cJ}{\mathcal{J}} 
\newcommand{\cX}{\mathfrak{X}} 
\newcommand{\cR}{\mathcal{R}} 
\newcommand{\cM}{\mathcal{M}} 
\newcommand{\cL}{\mathcal{L}}

\begin{document}

\title{Intrinsic Dynamics of Manifolds: 
Quantum Paths, Holonomy, and Trajectory Localization}  

\author{Mikhail Karasev\thanks{This research was supported 
by the RFBR (02-01-00952)
and by the INTAS (00-257).
The preprint version of the paper is in the arXiv: 
math.QA/0308163.}\\
\\
\small Moscow Institute of Electronics and Mathematics\\
\small Moscow 109028, Russia\\
\small karasev@miem.edu.ru}

\date{}

\maketitle

\begin{abstract}
We consider a dynamic generalization of the classical
``kinematic'' notion of affine connection providing a
correspondence between paths in the manifold and diffeomorphisms
of the manifold.  The linearization of these
path-diffeomorphisms coincides with the parallel translations
via the connection.  In this dynamic geometry one can translate
nonanalytic functions and distributions rather than tangent
vectors.  We describe the dynamic holonomy and the dynamic
curvature.

On the symplectic or quantum level this construction makes up
the symplectic or quantum paths, as well quantum  connection,
quantum curvature and quantum holonomy.

The construction of path-diffeomorphisms, 
being applied to trajectories of a given dynamical system, 
produces a dynamic localization of the system. 
In the symplectic (quantum) case, 
this dynamic localization provides 
a coherent-type representation of the quantum flow.
\end{abstract}

\maketitle

\section{Introduction}

Kinematics of manifolds is a parallel translation of tangent  
vectors along paths. The notion of affine connection
provides the mathematical description of this mechanical
concept. Usual constructions of the connection theory
(exponential mappings, geodesics, curvature, holonomy, etc.) 
are objects of the infinitesimal geometry. 
There are no potential forces in this geometry and thus 
no opportunities to translate nonanalytic functions 
or distributions along paths. 
Each path in the manifold knows only about 
the germ of the geometry but not 
about the neighborhood geometry of the manifold.
There is no place for the intrinsic wave concept 
in this kinematic approach. 
That is why one would like to ``integrate'' the usual connection
ideology to a more substantial framework, which can be called
a {\it dynamic geometry}. 
In particular, this generalization is necessary for solving the
quantization problem~\cite{K}. 

Dynamics is a translation of functions or distributions. 
The intrinsic dynamics of a manifold, 
which we have in mind, is generated by certain internal
non-autonomous vector field whose ``time'' variable ranges 
over the same manifold. 
In the symplectic case this internal vector field 
is Hamiltonian; the corresponding Hamilton
function determines the Ether structure~\cite{K}. 

In general, the internal field is a family of linear
maps acting in pairs of spaces tangent to the manifold.
This  family can also be considered 
as a one-form with values in vector fields, that is, in the Lie
algebra of first-order differential operators on the manifold. 
The basic condition for this form is the equation of zero
curvature.  

After restriction to the diagonal, where the ``time'' variable
coincides with the ``space'' variable, this zero curvature
equation occurs to be the Cartan structural equation (and its
solvability is guaranteed by the $1$-st Bianchi's identity).
Thus the idea of the internal field is a natural extension 
of Cartan's ``moving frame'' approach.
The main goal and advantage of such an extension is to get 
all charms and conveniences of zero curvature over manifolds
with a nonflat affine connection. 

Because of the zero curvature condition, 
the internal translation of functions (distributions,
curves, surfaces, etc.) 
along any path in the manifold does not depend on the shape of
the path. 
In this way we obtain a big family of diffeomorphisms
generating the first level of intrinsic dynamics of the
manifold. 

In the symplectic case, these diffeomorphisms are symplectic,
and among them there are Ether translations and
reflections~\cite{K}. They are applied to calculation of 
quantum geometric objects and determine, for instance,
the membrane representation of functions and the geometric phase
product~\cite{K1}.

In the given paper we consider a second level of intrinsic
dynamics. 
It is generated by the internal vector field with the
factor~$\frac12$. In presence of this factor one looses the
zero curvature condition:
the translation along a path now depends on the shape
of the path and determines a  transformation of the
manifold. We obtain a realization of the path groupoid
by diffeomorphisms.

Each path-diffeomorphism maps the origin of the path to its
end. The linearization of the path-diffeomorphism 
along the path itself coincides with the parallel translation by
means of a connection in the tangent bundle.
This is a point. 
We see that the path-diffeomorphism 
can be considered as a dynamic generalization of
the usual kinematic parallel translation of vectors.
In the symplectic framework, 
this generalization provides a notion of symplectic path 
and quantum path.

Applying the construction of path-diffeomorphisms and 
considering closed paths, at each given point 
we obtain a {\it dynamic holonomy\/} realized by 
diffeomorphisms with the given fixed point. Let us stress that
this is not the habitual holonomy in the tangent bundle.
Only the linearization at the fixed point gives us the usual
kinematic holonomy.  

The dynamic generalization of the notion of holonomy
opens the opportunity to obtain 
(in the symplectic case) the construction of quantum holonomy,
say, following the quantization approach~\cite{K,KM},
which allows one to quantize not only functions, but also
symplectic transformations.

In the standard way one can calculate the Lie algebra of the
dynamic holonomy group. This Lie algebra is related to certain
{\it dynamic curvature}. The germ of the dynamic
curvature at the fixed point is given by the ordinary
kinematic curvature (and by the symplectic form if we are
inside the symplectic framework). 

At the next stage we consider a certain external dynamical 
(Hamiltonian) system on the manifold. 
The trajectories of this system are paths, 
and to them one can assign path-diffeomorphisms
(symplectomorphisms) of the manifold. 
Each diffeomorphism can be applied to any function or 
to any dynamical system on the manifold. 
In particular, if it is applied to the given
system itself, then one obtains a new system 
that has the equilibrium point at the origin of the trajectory. 

In a sense, this is a procedure of localization of the
dynamical system at the given trajectory. 
That is why we use the term ``localization'' 
in the title of this paper.
But we would like to stress that 
the new transformed system not only inherits the germ 
(say, the first variation) of the original system
at the chosen trajectory, 
but provides the total global information about the original
system.  
For this type of transformations, we use
the name {\it translocations}. 

The translocation is determined by a chosen trajectory of a
dynamical system and by the intrinsic dynamics of the whole
manifold. This is a representation of the system via the
intrinsic travelling stream generated by the internal vector 
field on the manifold.

In the symplectic case the translocation reduces a given
Hamiltonian system to a new one with a Hamilton function 
approximately quadratic near the given point. The corresponding
quantum transformation is well known in the Euclidean
space and is extensively used in the oscillatory approximation
for wave and quantum equations(see, for instance,~\cite{Bab}). 
In the case of general symplectic manifolds, the existence of such
transformations was in question.
The concept of Ether structure~\cite{K} allows one to answer 
this question. 

In this paper we consider both the classical and quantum 
aspects of these problems and develop the {\it dynamic
geometry\/} framework following the ideology 
of~\cite{K}. 
The symplectic case is of our main interest here, and therefore
the basic text of the paper deals with this case. 
The Appendix contains a parallel description of the general
situation (manifolds with affine connection).

\section{Ether dynamics}

Let $\cX$ be a manifold with a symplectic form~$\omega$, and
with a torsion free symplectic connection~$\Gamma$.
By $\{f,g\}=\nabla f\cdot\Psi\cdot \nabla g$ 
we denote the Poisson brackets related to~$\omega$, and
by~$\nabla$ the covariant derivative with respect to~$\Gamma$.

The {\it Ether structure\/} on $\cX$ is 
given by a certain intrinsic Hamiltonian~$\cH$ which is a
$1$-form on~$\cX$ with values in smooth real functions on~$\cX$.
In local coordinates:
$$
\cH_x(z)=\sum^{2n}_{j=1}\cH_x(z)_j\,dx^j,\qquad 2n=\dim\cX,
\quad x,z\in\cX.
$$
The following conditions should hold:

(a) the zero curvature condition
\begin{equation}
\pa\cH+\frac12\{\cH\wedgeco\cH\}=0
\tag{2.1}
\end{equation}
(where $\pa=d_x$ is the differential by~$x$ 
and the Poisson brackets are taken by~$z$);

(b) the boundary conditions
\begin{equation}
\cH\Big|_{\diag}=0,\qquad
\nabla\cH\Big|_{\diag}=2\omega,\qquad
\nabla\nabla\cH\Big|_{\diag}=0
\tag{2.2}
\end{equation}
(where $\nabla=\nabla_z$ and $\diag=\{z=x\}$);

(c) the skew-symmetry condition
\begin{equation}
\cH_x\big(s_x(z)\big)=-\cH_x(z),
\tag{2.3}
\end{equation}
where $s=s_x(z)$ is the solution of the Ether Hamiltonian system 
\begin{equation}
\pa s=\nabla\cH(s)\Psi(s),\qquad s\Big|_{x=z}=z.
\tag{2.4}
\end{equation}

In the last system $\pa=d_x$ and the variable $x\in\cX$ plays the
role of time. 

\begin{theorem}
{\rm(i)} For any symplectic connection on a symplectic
manifold~$\cX$, the zero curvature equation {\rm(2.1)}
with conditions {\rm(2.2)}, {\rm(2.3)} is solvable at
least in a neighborhood of the diagonal $\diag=\{z=x\}$.

{\rm(ii)} The mappings $s_x$ defined by {\rm(2.4)}
introduce a reflective structure to~$\cX${\rm:}

-- $s_x$ are symplectic,

-- $s^2_x=\id$,

-- $s_x$ has the isolated fixed point~$x$.

{\rm(iii)} By a given reflective structure 
on the simply connected symplectic manifold~$\cX$, 
one reconstructs the corresponding symplectic
connection and the Ether Hamiltonian as follows{\rm:}
\begin{align*}
\Gamma^j_{kl}(x)&=-\frac12\frac{\pa^2 s^j_x(z)}{\pa z^k\pa z^l}
\bigg|_{z=x},\\
\cH_x(z)_k&=\int^z_x\frac{\pa s^j_x}{\pa x^k}
\big(s_x(z)\big)\omega_{jm}(z)\,dz^m,
\end{align*}
where the integral is taken along any path connecting~$x$ with~$z$.

{\rm(iv)} In addition to the boundary conditions {\rm(2.2)},
the third covariant derivative of $\cH_x(z)$ by~$z$ 
on the diagonal $z=x$ is related to the curvature tensor~$R$ of the
connection~$\Gamma${\rm:} 
$$
\nabla_j\nabla_l\nabla_s \cH_k\Big|_{\diag}=2\omega_{sm}R^m_{lkj}.
$$
All higher covariant derivatives of $\cH$ on the diagonal are also
explicitly evaluated in terms of~$\omega$ and~$R$.
\end{theorem}

Note that the zero curvature condition (2.1) guarantees the
solvability of (2.4). Besides of (2.4), one can consider more
general Hamiltonian system
\begin{equation}
\frac{d}{dt}G=\dot y(t)\nabla \cH_{y(t)}
\big(G\big)\Psi\big(G\big),\qquad G\bigg|_{t=0}=x.
\tag{2.5}
\end{equation}
Here $\{y(t)\mid 0\leq t\leq 1\}$ is a smooth path connecting 
a point $y_0=y(t)|_{t=0}$ with a point $y_1=y(t)|_{t=1}$. 
By $\dot y(t)$ we denote the velocity vector at the point $y(t)$
of the path. 

The zero curvature condition (2.1) implies that the solution
$G=G(t)$ of (2.5) at $t=1$ does not depend on 
the shape of the path $\{y(t)\}$, that is, we can denote:
$$
G\big|_{t=1}=g_{y_1,y_0}(x).
$$
Thus one obtains a family of symplectic transformations
$g_{y_1,y_0}$ of the manifold~$\cX$. We call them~\cite{K}
the {\it Ether translations}.
They are related to the reflections as follows.

\begin{theorem}
$\quad g_{y,x}=s_y s_x$, $\quad g_{y,x}(x)=s_y(x)$.
\end{theorem}

\section{Factor $\frac12$}

The symplectic connection $\Gamma$ has its geodesics and
exponential mappings $\exp_x$, defined in the standard way.
The Ether Hamiltonian~$\cH$, which provides an ``integral''
approach to the connection theory on~$\cX$, generates its own
geodesics and exponential mappings by means of the Hamiltonian
translations. Namely, let us consider the Hamiltonian system~\cite{K}:
\begin{equation}
\frac{d}{dt}E=\frac12 v \nabla \cH_x(E) \Psi(E),\qquad
E\bigg|_{t=0}=x.
\tag{3.1}
\end{equation}
Here $v\in T_x\cX$ is a fixed velocity vector. Denote the
solution of (3.1) by
$$
E(t)=\Exp_x(vt).
$$
This trajectory we call the {\it Ether geodesics\/} through~$x$,
and $\Exp_x$ is called the Ether exponential mapping. 
This mapping, in general, is different from $\exp_x$, and the
Ether geodesics do not coincide with the $\Gamma$-geodesics.

The following relations hold:
\begin{align*}
s_x\big(\Exp_x(v)\big)&=\Exp_x(-v),\\
\cH_x\big(\Exp_x(v)\big)&=-\cH_x\big(\Exp_x(-v)\big).
\end{align*}

Let us now pay attention to the factor~$\frac12$ on the right 
of (3.1). We can consider more generic systems with such a
factor: 
\begin{equation}
\frac{d}{dt}Y=\frac12 \dot y(t)\nabla\cH_{y(t)}(Y)\Psi(Y),\qquad
Y\bigg|_{t=0}=x.
\tag{3.2}
\end{equation}
Here $\{y(t)\}$ is a given smooth path in~$\cH$. 
The solution of (3.2) we denote by $Y=Y^t(x)$;
thus $Y^t$ is the translation along the trajectories of~(3.2). 

Of course, the mapping $Y^t$, determined by the Hamiltonian
system, is symplectic.
But in the presence of the factor~$\frac12$ in (3.2), 
we loose the zero curvature property, and so, 
in contrast to~(2.5),
the mapping~$Y^t$ does depend on the shape of the path 
$\{y(t)\}$. 
Let us formulate the key observation about this mapping. 

\begin{theorem}
The following identities hold{\rm:}
$$
Y^t(y(0))=y(t),\qquad dY^t(y(0))=V^t.
$$
Here $V^t$ is the parallel translation 
$T_{y(0)}\cX\to T_{y(t)}\cX$ by means of 
the connection~$\Gamma$ along the path 
$\{y(\tau)\mid 0\leq \tau\leq t\}$.
\end{theorem}

{\it Proof}.
In view of the boundary conditions~(2.2),
we have
$$
\frac12 \nabla\cH_y(Y)\Psi(Y)=I+\langle \Lambda,\delta Y\rangle
+O(\delta Y^2).
$$
Here $\delta Y=Y-y$ in some local coordinates, and 
$\Lambda=\omega(\pa \Psi+\Gamma\Psi)$.
Taking the symplecticity of the connection~$\Gamma$ into
account,  the system (3.2) can be written as 
\begin{equation}
\frac{d}{dt}\delta Y+\dot y\Gamma(y)\,\delta Y
+O(\delta Y^2)=0.
\tag{3.3}
\end{equation}

If the initial data for $Y$ coincides with the initial data
for~$y$, then $\delta Y|_{t=0}=0$, and from (3.3) 
we have $\delta Y=0$ for any~$t$. 
Thus the first identity of the theorem is true.

From (3.3) it also follows that the differential 
$V^t=d Y^t(y(0))$ obeys the equation
$$
\frac{d}{dt}V^t+\dot y(t)\Gamma(y(t))\, V^t=0.
$$
This means that $V^t$ is the parallel translation along the path
$\{y(\tau)\mid 0\leq \tau\leq t\}$ by means of the
connection~$\Gamma$.  
The theorem is proved.

This theorem explains why we take the factor~$\frac12$ on the
right-hand side of~(3.2). Any different choice destroys the
statement of the theorem.

\section{Symplectic paths}

Let us introduce more convenient notation. 
A path in~$\cX$ we will denote by $\sigma=\{y(t)\}$, 
and the translation along trajectories of (3.2)
denote by $[\sigma]=Y^t$. 
So, $[\sigma]\big(y(0)\big)=Y$ is the solution of (3.2).

\begin{definition}
The symplectic transformation $[\sigma]:\,\cX\to\cX$ we call 
the {\it symplectic path\/} corresponding to the classical
path~$\sigma\subset \cX$. 
\end{definition}

Consider the induced mapping 
$[\sigma]^*:\, C^\infty(\cX)\to C^\infty(\cX)$. From equation
(3.2) we have the following explicit formula:
\begin{equation}
[\sigma]^*=\underset{\to}{\Exp}\bigg(\frac12\int_\sigma \ad(\cH)\bigg).
\tag{4.1}
\end{equation}
Here $\ad(f)$ denotes the Hamiltonian vector field on~$\cX$ 
related to a function~$f$, and $\underset{\to}{\Exp}$ 
is the multiplicative exponential 
ordered from left to right:
\begin{equation}
\underset{\to}{\Exp}\bigg(\int^t_0 a\bigg)
=\lim_{N\to\infty}e^{\Delta t\, a(t_1)}\dots e^{\Delta t\, a(t_N)},
\tag{4.2}
\end{equation}
where $\Delta t =t/N$ and $t_j=j\Delta t$.
Formula (4.1) represents the parallel translation 
along the path~$\sigma$ with respect to the connection
\begin{equation}
{\boldsymbol\nabla}^0=\pa+\frac12 \ad(\cH)
\tag{4.3}
\end{equation}
acting in the trivial bundle over~$\cX$ 
with fibers $C^\infty(\cX)$.

\begin{theorem}
{\rm(i)}
The correspondence $\sigma\to[\sigma]$ between classical paths
and symplectic paths is given by formula~{\rm(4.1)}. 
This is a realization of the path groupoid in the group of
symplectic transformations of~$\cX${\rm:}
if the product of two paths $\sigma_2\circ\sigma_1$ exists, then
$$
[\sigma_2\circ\sigma_1]=[\sigma_2]\circ[\sigma_1].
$$

{\rm(ii)}
The differential $d[\sigma]$ of the symplectic path 
at the origin of the path coincides
with the parallel translation along the path~$\sigma$
by means of the connection~$\Gamma$.

{\rm(iii)}
Symplectic paths commute with reflections, i.e., 
if $\sigma_{y,x}$ is a path connecting~$x$ with~$y$, then
\begin{equation}
[\sigma_{y,x}]\circ s_x=s_y\circ[\sigma_{y,x}].
\tag{4.4}
\end{equation}
In particular, the mapping $[\sigma_{y,x}]$ is affine along the
path $\sigma_{y,x}$, that is, 
it translates the connection coefficients $\Gamma(x)$ to
$\Gamma(y)$. 
\end{theorem}

\section{Symplectic loops and Ether curvature}

Now, for any point $x\in\cX$, we can consider the group 
of closed paths~$\sigma$  (parametrized loops) starting
from~$x$. The corresponding symplectic loops~$[\sigma]$
form a subgroup in the whole group of symplectic transformations
of the manifold~$\cX$ having~$x$ as a fixed point.
This subgroup we call a {\it dynamic holonomy group\/} 
and denote by~$\cJ_x$. 

The differentials $d[\sigma](x)$ of symplectic loops 
at the fixed point~$x$ are symplectic linear transformations 
of $T_x\cX$. They form a subgroup in the whole group of linear 
transformations of $T_x\cX$.
This subgroup coincides with the usual holonomy
group of the connection~$\Gamma$ at the point~$x$.
The latter group we shall call {\it kinematic\/}
in order to distinguish it from the dynamic holonomy group. 

Obviously, the dynamic holonomy groups $\cJ_x$ and $\cJ_y$
corresponding to different points~$x$ and~$y$ are gauge
equivalent (conjugate to each other):
$$
\cJ_x=[\sigma_{y,x}]^{-1}\circ\cJ_y\circ[\sigma_{y,x}],
$$
where $\sigma_{y,x}$ is a path connecting~$x$ with~$y$.

Let us mark a certain point $0\in\cX$. 
Consider small membranes $\Sigma\subset\cX$ whose boundaries are
parametrized loops $\pa\Sigma\in\cJ_0$. 

Note that the multiplicative exponential~(4.2) can be
represented (using the continual version of the
Campbell--Hausdorff formula) as follows:
$$
\underset{\to}{\Exp}\bigg(\int_{\pa \Sigma} a\bigg)
=\exp\bigg(\int_{\Sigma} b+o(\Sigma)\bigg).
$$
Here $b=da+\frac12[a\wedge a]$ is the curvature of~$a$, 
and the summands $o(\Sigma)$ are of higher degree
with respect to the area of~$\Sigma$.

Applying this fact to (4.1), we see that 
$$
[\pa\Sigma]^{-1*}=\exp\bigg(-\frac12\int_\Sigma\ad
\Big(\pa\cH+\frac14\{\cH\wedgeco\cH\}\Big)+o(\Sigma)\bigg).
$$
The sign minus in the exponent is due to the additional
inversion of the map $[\pa\Sigma]$; we need it to have the
group homomorphism $\pa\Sigma\to [\pa\Sigma]^{-1*}$.

Taking into account the zero curvature condition~(2.1), 
we derive 
\begin{equation}
[\pa\Sigma]^{-1*}=\exp\bigg(\frac18\int_\Sigma\ad
(\{\cH\wedgeco\cH\})+o(\Sigma)\bigg).
\tag{5.1}
\end{equation}
These are operators in $C^\infty(\cX)$ representing the 
dynamic holonomy group~$\cJ_0$. 

Near the unity element $0\in \cJ_0$ any loop $\sigma\in\cJ_0$ can
be considered as a plane one.
The differential of the family of operators (5.1) at the unity
point of the group~$\cJ_0$ is the operator-valued $2$-form
\begin{equation}
d[\sigma]^{-1*}\bigg|_{\sigma=0}
=\frac18\ad(\{\cH_0\wedgeco\cH_0\}).
\tag{5.2}
\end{equation} 
Here $\cH_0$ is the Ether Hamiltonian at the marked point
$0\in\cX$.  

Let us fix certain local coordinates near the marked point
$0\in\cX$. 
For any small $\sigma\in\cJ_0$ denote by $\sigma_{jk}$ the area
enclosed by the projection of $\sigma$ onto $(j,k)$-coordinate
plane (with the orientation given by $dx^k\wedge dx^j$). 
Denote by $\cH_{0j}$ the components of the Ether Hamiltonian
$\cH_0$ with respect to the chosen local coordinates, and  
introduce the following functions: 
\begin{equation}
\cR_{jk}\od\frac14\{\cH_{0k},\cH_{0j}\}.
\tag{5.3}
\end{equation}

Note that from the zero curvature equation (2.1) and from 
definition (4.3) of the connection ${\boldsymbol\nabla}^0$
we know that the Hamiltonian fields of the functions 
$\cR_{jk}$ 
represent the curvature of ${\boldsymbol\nabla}^0$, i.e.,
$$
[{\boldsymbol\nabla}^0_j,{\boldsymbol\nabla}^0_k]=\ad(\cR_{jk}).
$$

From (5.2) one has the formula
\begin{equation}
\frac{\pa}{\pa\sigma_{jk}}[\sigma]^{-1*}\bigg|_{\sigma=0}=\ad(\cR_{jk}).
\tag{5.4}
\end{equation}
Thus we obtain the following statement.

\begin{theorem}
If $\cX$ is simply connected, 
then the Lie algebra of the dynamic holonomy group~$\cJ_0$ 
at the marked point $0\in\cX$ is generated by
Hamiltonian vector fields related to functions~$\cR_{jk}$ 
that are defined in~{\rm(5.3)}.

The marked point is a stationary point of the 
Hamiltonians $\cR_{jk}$, 
namely:
$$
\cR_{jk}(0)=\omega_{jk}(0),\qquad \nabla \cR_{jk}(0)=0.
$$
Their second derivatives  are given by the curvature
tensor~$R$ of the connection~$\Gamma$ as follows{\rm:} 
$$
\nabla^2_{lm} \cR_{jk}(0)=2\omega_{ls}(0) R^s_{mjk}(0).
$$
\end{theorem}

We note that the curvature matrices 
$R_{jk}(0)=(\!(R^s_{mjk}(0)\,)\!)$ 
represent operators in $T_0\cX$ which are skew-symmetric with respect to
the symplectic inner product generated by $\omega(0)$, that is,
they belong to the Lie algebra $\spp(T_0\cX)$.

\begin{corollary}
The Lie algebra of the symplectic holonomy 
group of~$\cX$ at the marked point~$0$
is generated by components of the curvature
form $R_{jk}(0)\in \spp(T_0\cX)$.
\end{corollary}

This statement is very common 
for the ``kinematic'' connection theory:
it is a particular case of the general Ambrose--Singer theorem
and is analogous to of E.~Cartan's theorem about the Riemannian
holonomy. 
In Theorem~5.1 we have its generalization 
to the dynamic holonomy group.

We shall call the form 
$\cR=\frac18\{\cH\wedgeco\cH\}$ 
staying in formula~(5.1) the {\it Ether curvature form}.

\section{Translocations}

Let us consider now, on the symplectic manifold~$\cX$, 
an exterior Hamiltonian system related to a smooth function~$H$:
\begin{equation}
\frac{d}{dt}X=\nabla H(X)\,\Psi(X),\qquad X\bigg|_{t=0}=x.
\tag{6.1}
\end{equation}
Denote the solution by $X=X^t(x)$; thus $X^t$ is the Hamilton
flow generated by~$H$.

Let us fix a point $y\in\cX$ and a time $t\in\bR$.
We shall use the same notation 
both for the point $X^t(y)$ of the trajectory 
and for the whole segment of the trajectory:
\begin{equation}
X^t(y)\sim\{X^\tau(y)\mid 0\leq \tau\leq t\}.
\tag{6.2}
\end{equation}
The symplectic path corresponding to the path~(6.2)  
is denoted by $[X^t(y)]$
and will be called a {\it symplectic trajectory}. 
For each fixed~$t$,  
it is a symplectic transformation of~$\cX$.

We can define a new Hamilton function
\begin{equation}
H^t_y\od [X^t(y)]^*\Big(H-\frac12
\dot X^t(y)\,\cH_{X^t(y)}\Big)-H(y).
\tag{6.3}
\end{equation}
Related Hamiltonian system is 
\begin{equation}
\frac{d}{dt}Z=\nabla H^t_y(Z)\,\Psi(Z),\qquad Z\bigg|_{t=0}=x.
\tag{6.4}
\end{equation}
Denote the solution of this system by $Z=Z^t_y(x)$;
thus $Z^t_y$ is the translation along trajectories of~(6.4).

\begin{theorem}
{\rm(i)} The translation $X^t$ along trajectories of the
original Hamiltonian system {\rm(6.1)} 
and the translation~$Z^t_y$ along trajectories of the
transformed Hamiltonian system {\rm(6.4)} 
relate to each other as follows{\rm:}
\begin{equation}
X^t=[X^t(y)]\circ Z^t_y.
\tag{6.5}
\end{equation}
Here $y\in\cX$ is an arbitrary point, and $[X^t(y)]$ is the
symplectic trajectory corresponding to the segment
{\rm(}path{\rm)} {\rm(6.2)}.

{\rm(ii)} The point~$y$ is a stationary point of
the Hamiltonian $H^t_y$ of system {\rm(6.4)}, 
namely{\rm:} 
$H^t_y(y)=0$, $(\nabla H^t_y)(y)=0$.
The second derivative matrix at~$y$ is
\begin{equation}
\nabla^2 H^t_y(z)\Big|_{z=y}
= V^{t*}_y\cdot\nabla^2 H\big(X^t(y)\big)\cdot V^t_y.
\tag{6.6}
\end{equation}
Here $\nabla^2 H$ is the tensor of second covariant derivatives
of~$H${\rm:} 
\begin{equation}
\nabla^2_{jk} H \od D^2_{jk}H-D_mH\cdot \Gamma^m_{jk},
\tag{6.7}
\end{equation}
and $V^t_y$ is the parallel translation along the path {\rm(6.2)}
by means of the connection~$\Gamma$.

{\rm(iii)} The symplectic mapping $Z^t_y$ in {\rm(6.5)} 
has the fixed point~$y$. The symplectic mapping $[X^t(y)]$
transports~$y$ to $X^t(y)$.

The differential of the solution to system {\rm(6.1)} at the
point~$y$ is given by 
\begin{equation}
dX^t(y)=V^t_y\circ W^t_y.
\tag{6.8}
\end{equation}
Here $W=W^t_y$ is the solution of the linear system over $T_y\cX$:
\begin{equation}
\frac{d}{dt}W=M^t_y W,\qquad W\bigg|_{t=0}=I,
\tag{6.9}
\end{equation}
where
\begin{equation}
M^t_y\od -\Psi(y)\cdot V^{t*}_{y}\cdot \nabla^2 H(X^t(y))\cdot V^t_y.
\tag{6.10}
\end{equation}
\end{theorem}

The transformation from the system (6.1) to the system (6.4)
will be called a {\it translocation\/} to~$y$.
This is a dynamic analog of the usual kinematic 
parallel translation. 

The translocation localizes the whole system 
near the given trajectory starting at~$y$. 
In particular, the first variation 
of the original system along the trajectory 
is reduced to the linear equation (6.9)
whose matrix is determined by the second covariant derivatives 
of~$H$ via (6.10) 
(about transformations of the first variation system
via symplectic connections see a detailed investigation
in~\cite{MRR}). 

\begin{corollary}
Let the Hamilton function~$H$ be covariantly quadratic 
along the trajectory {\rm(6.2)}, i.e., 
\begin{equation}
\nabla_{\nabla H\Psi}(\nabla^2 H)=0.
\tag{6.11}
\end{equation}
at all points of the trajectory {\rm(6.2)}.
Then the solution of the first variation system along this
trajectory is given by
\begin{equation}
dX^t(y)=V^t_y\circ\exp\big(-t\Psi(y)\nabla^2 H(y)\big).
\tag{6.12}
\end{equation}
If the trajectory {\rm(6.2)} is periodic, then its monodromy
matrix at~$y$ is factorized to the product of the geometric
monodromy {\rm(}the holonomy of the connection{\rm)}
and the dynamical monodromy 
{\rm(}with the generator $-\Psi(y)\nabla^2 H(y)${\rm)}.
\end{corollary}

\section{Quantum paths and quantum curvature}

We assume that there is a representation of a space of functions
on $\cX$ in a Hilbert space given by the integral
\begin{equation}
f\to\hat{f}=\frac1{(2\pi\hb)^n}\int_{\cX}f\bS\,dm.
\tag{7.1}
\end{equation}
Here $\hb>0$, $2n=\dim\cX$, $dm$ is a measure on $\cX$, 
and $\bS=\{\bS_x\mid x\in\cX\}$ is a family of operators in the
Hilbert space obeying the Schr\"odinger type dynamics equation
over~$\cX$:  
\begin{equation}
i\hb \pa \bS=\wh{\cH}^\hb\bS,
\tag{7.2}
\end{equation}
where $\cH^\hb$ is a quantum intrinsic Hamiltonian 
such that the quantum zero curvature equation holds:
\begin{equation}
\pa \wh{\cH}^\hb+\frac{i}{2\hb}
[\wh{\cH}^\hb\wedgeco\wh{\cH}^\hb]=0.
\tag{7.3}
\end{equation}

As in the work~\cite{K}, 
in addition to (7.2), (7.3), 
we assume that $\wh{\cH}^\hb$ is self-adjoint 
and all $\bS_x$ are self-adjoint and almost unitary operators: 
\begin{equation}
(\wh{\cH}^\hb)^*=\wh{\cH}^\hb,\qquad
\bS^*_x=\bS_x,\qquad 
\bS^2_x=\mu^2\cdot I,
\tag{7.4}
\end{equation}
where $\mu=\const$.  
Note that the first and third conditions in (7.4) are not necessary, 
but we fix them to simplify all formulas. 
Under these conditions the term of order~$\hbar$
in the expansion of the quantum Hamiltonina $\cH^\hbar$ is absent;
namely, 
we have 
$$
\cH^\hb=\cH+O(\hb^2),
$$
where $\cH$ is a classical intrinsic Hamiltonian
from Sect.~2.

In this situation the quantum versions of reflections $s_x$ and 
Ether exponential mappings $\Exp_x$ were describe in~\cite{K},
and the semiclassical asymptotic formulas were derived.

Now the first question arises: 
what is the quantum version of path-symplectomorphisms? 

The answer is:
\begin{equation}
\wh{U}_\sigma =\underset{\leftarrow}{\Exp}
\bigg(-\frac{i}{2\hb}\int_{\sigma}\wh{\cH}^\hb\bigg).
\tag{7.5}
\end{equation}
This unitary operator acting in the Hilbert space will be called
the {\it quantum path\/} corresponding to the classical path 
$\sigma\subset\cX$.
It is generated by the {\it quantum connection\/}
$$
{\boldsymbol\nabla}^\hb=\pa+\frac{i}{2\hb}\wh{\cH}^\hb.
$$

Of course, the correspondence $\sigma\to\wh{U}_\sigma$ 
is a representation of the path groupoid:
$$
\wh{U}_{\sigma_2}\cdot\wh{U}_{\sigma_1}
=\wh{U}_{\sigma_2\circ\sigma_1},
\qquad\text{or}\qquad
U_{\sigma_2}*U_{\sigma_1}=U_{\sigma_2\circ \sigma_1}.
$$
Here $*$ is the noncommutative product in the function algebra
over $\cX$ which corresponds to the operator representation
(7.1): 
$$\hat{f}\hat{g}=\wh{f*g}
$$
(see details and references in~\cite{K}).

Formula (7.5) is a Schr\"odinger type representation of the
path~$\sigma$,  
but one can define the Heisenberg type representation as well:
$$
\wh{U}_\sigma^{-1}\cdot \hat{f}\cdot \wh{U}_\sigma
=\wh{\hat{\sigma}^*f}.
$$
The mapping $\hat{\sigma}^*$ acts in the noncommutative algebra
of functions over $\cX$ by the formula 
\begin{equation}
\hat{\sigma}^*=\underset{\to}{\Exp}
\bigg(\frac{i}{2\hb}\int_{\sigma}\ad_{\cH^\hb}\bigg),
\qquad
\ad_{\cH^\hb}(f)\od\cH^\hb*f-f*\cH^\hb,
\tag{7.6}
\end{equation}
and possesses the $\hb$-expansion:
\begin{equation}
\hat{\sigma}^*=[\sigma]^*\big(I+O(\hb^2)\big),
\tag{7.7}
\end{equation}
where $[\sigma]$ is the classical path-symplectomorphism (4.1).

The definition (7.5) and formulas (7.6), (7.7)  
allow one to determine the {\it quantum holonomy} groups 
$\wh{\cJ}_x$ (groups of unitary operators) and to calculate 
their Lie algebras. The quantum version of Theorem~5.1 is the
following one.

\begin{theorem}
Let $\cX$ be simply connected. 
For any point $x\in\cX$ the Lie algebra of the quantum holonomy
group $\wh{\cJ}_x$ is generated by the quantum curvature operators
$$
\wh{\cR}^\hb_{xjk}
=\frac{i}{4\hb}\big[\,\wh{\cH}^\hb_{xk},\wh{\cH}^\hb_{xj}\,\big],
$$
where $\wh{\cH}^\hb_{xk}$ are components 
of the quantum intrinsic Hamiltonian. 
\end{theorem}

The {\it quantum curvature form\/} is naturally defined as 
$$
\wh{\cR}^\hb
=\frac{i}{8\hb}\big[\wh{\cH}^\hb\wedgeco\wh{\cH}^\hb\big].
$$
The quantum analog of (5.1) is the following:
\begin{equation}
\big(\wh{U}_{\pa\Sigma}\big)^{-1}
=\exp\bigg(-\frac{i}{\hbar}\int_{\Sigma}\wh{\cR}^\hb+o(\Sigma)\bigg).
\tag{7.8}
\end{equation}

Following \cite{K}, we can derive the simplest semiclassical
asymptotics of quantum paths, assuming, for simplicity,  
that the form~$\omega$ is exact. 

Let $\sigma \subset \cX$ be a path, 
and let for any $x\in\cX$ the symplectic mapping 
$s_x\circ[\sigma]$ have a unique fixed point~$\tilde{x}$
(smoothly depending on~$x$). 
Consider the membrane $\Sigma_\sigma(x)$ 
whose boundary consists of four pieces: 
a path~$c$ from $\tilde{x}$ to the beginning point of~$\sigma$, 
the path~$\sigma$,
the path $[\sigma](c)$, 
and the Ether geodesics connecting $[\sigma](\tilde{x})$ 
with~$\tilde{x}$ through the midpoint~$x$.

\begin{theorem}
{\rm(i)}
Under the above assumptions, the symbol $U_\sigma$ of the
quantum path {\rm(7.5)} has the asymptotics
\begin{equation}
U_\sigma(x)=\exp\bigg\{\frac i\hb
\int_{\Sigma_\sigma(x)}\omega\bigg\}
\varphi_\sigma(x)+O(\hb),
\tag{7.9}
\end{equation}
where
\begin{equation}
\varphi_\sigma(x)=2^n\cdot
\det \big(I-D(s_x\circ[\sigma])(\tilde{x})\big)^{-1/2}.
\tag{7.10}
\end{equation}

{\rm(ii)}
The membrane $\Sigma_\sigma(x)$ in {\rm(7.9)} 
can be chosen as a union of two pieces. 
The first piece $\Sigma_\sigma$
is independent of~$x$; 
its boundary is formed by the path $\sigma$ 
and the Ether geodesic connecting the end-points 
of~$\sigma$.
Let us denote by $0_\sigma$ the mid-point 
of this Ether geodesic.
The second piece is bounded by two paths~$c$ 
and~$[\sigma](c)$,
connecting the points $\tilde{x}$ 
and $[\sigma](\tilde{x})$
with the end-points of $\sigma$, 
and by two Ether geodesics
through the mid-points~$x$ and $0_\sigma$.
The WKB-phase in {\rm(7.9)} 
is represented as the sum 
\begin{equation}
\int_{\Sigma_\sigma(x)}\omega
=\int_{\Sigma_\sigma}\omega
+\Phi_{0_\sigma}^{[\sigma]}(x),
\tag{7.11}
\end{equation}
where $\Phi^{[\sigma]}$ 
is the normalized generating function 
of $[\sigma]$ in the sense of~{\rm\cite{K1}}
{\rm(}see formula {\rm(3.3))}.
\end{theorem}

\begin{remark}
For any {\rm(}closed or nonclosed{\rm)} path~$\sigma$ the phase
\begin{equation}
\Phi^{[\sigma]}(x)=\int_{\Sigma_\sigma(x)} \omega
\tag{7.12}
\end{equation}
appearing in (7.9) is a generating function 
of the symplectic mapping~$[\sigma]$
in the sense of~\cite{K1}.
In particular, one has the relation
\begin{equation}
d\Phi^{[\sigma]}(x)=\cH_x([\sigma](\tilde{x})\,)=-\cH_x(\tilde{x}),
\tag{7.13}
\end{equation}
where $\tilde{x}$, as above, denotes the fixed point of the mapping 
$s_x\circ[\sigma]$.
One can rewrite relations (7.13) in the following form:
$$
[\sigma](\tilde{x})=l(x,d\Phi^{[\sigma]}(x)),\qquad
\tilde{x}=r(x,d\Phi^{[\sigma]}(x)).
$$
Here $l$ and $r$ are left and right mappings 
in the symplectic groupoid over $\cX$
(i.e., in a neighborhood of the zero section in $T^*\cX$), 
see~\cite{K,K1,KM}.

Also note that in view of (7.5) the function $U_\sigma$
satisfies the Cauchy problem for 
a Schr\"odinger-type evolution equation,
and so, the phase $\Phi^{[\sigma]}$ (7.9), (7.12) 
satisfies the Hamilton--Jacobi equation
\begin{equation}
\frac{\pa}{\pa\sigma''}\Phi^{[\sigma]}(x)
+\frac12 \cH_{\sigma''}\big(l(x,d\Phi^{[\sigma]}(x))\big)=0.
\tag{7.14}
\end{equation}
Here $\frac{\pa}{\pa\sigma''}$ means the differential 
with respect to variations of the final endpoint of the path 
$\sigma=\{\sigma'\to\sigma''\}$.
\end{remark}

\section{Quantum evolution via translocation}

Now let us answer the next question:
what is the quantum version of translocations? 

Let us consider an external Hamiltonian~$H$ on~$\cX$
and the corresponding quantum operator $\wh{H}$.
The {\it quantum translocation\/} is the unitary operator
\begin{equation}
\wh{U}^{t}_{y}
=\exp\bigg\{-\frac i\hb H(y)t \bigg\} \wh{U}_{X^t(y)},
\tag{8.1}
\end{equation}
where $\wh{U}_{X^t(y)}$ is the quantum trajectory 
(quantum path) corresponding to the segment of the classical
Hamiltonian trajectory (6.2).  

\begin{theorem}
The following permutation formula holds{\rm:}
\begin{equation}
\bigg(-i\hb\frac{\pa}{\pa t}+\wh{H}\bigg)\cdot \wh{U}^t_{y}
=\wh{U}^t_{y}\cdot
\bigg(-i\hb\frac{\pa}{\pa t}+\wh{\bH}^t_y\bigg),
\tag{8.2}
\end{equation}
where
$$
\wh{\bH}^t_y\od \wh{X^t(y)}^*
\Big(H-\frac12 \dot X^t(y)\cH^\hb_{X^t(y)} \Big)-H(y).
$$
The quantum translocated Hamiltonian $\bH^t_y$ differs from 
the classical translocated Hamiltonian $H^t_y$ {\rm(6.3)} 
by $O(\hb^2)${\rm:} 
\begin{equation}
\bH^t_y=H^t_y+O(\hb^2)\qquad \text{as}\quad \hb\to0.
\tag{8.3}
\end{equation}
\end{theorem}

As a corollary, we obtain the following representation 
of the quantum Schr\"odinger type flow:
\begin{equation}
\exp\bigg\{-\frac{it}{\hb}\wh{H}\bigg\}
=\wh{U}^t_y\cdot \underset{\leftarrow}{\Exp}
\bigg\{-\frac i\hb\int^t_0 {\wh{\bH}}^t_y\,dt\bigg\}.
\tag{8.4}
\end{equation}
From Theorems 6.1 and 7.2 and formulas (8.1), (8.3),
one can now easily derive the global semiclassical evolution of
any quantum state localized at the point~$y$.

Indeed, the operator $\wh{\bH}^t_y=\wh{H}^t_y+O(\hb^2)$ 
in (8.4), near the point~$y$, 
looks like an oscillator Hamiltonian perturbed by cubic and
higher-order terms. 
Thus if the Wigner function of the initial state is localized
near~$y$ (for example, if it is a Gaussian type function), 
then this localized behavior will not change after application
of the evolution operator 
$\Exp\{-\frac i\hb\int^t_0 \wh{\bH}^t_y\,dt\}$.
For example, in the Gaussian case, 
as it follows from the general analysis \cite{Bab,K2,Gui,MS,VPM},
the quadratic part of the initial Gaussian exponent will remain
quadratic and just transformed by the first variation system
related to the Hamiltonian $H^t_y$; 
in our case this is system (6.9) over $T_y\cX$.

At the last stage one can integrate over all points $y\in\cX$
and calculate the quantum evolution for generic Cauchy data
(as in~\cite{P}).

This is just the general scheme. 
Now let us present more details.
Let $\Pi^0_{y,\hb}$ be a family of functions 
(coherent states) over~$\cX$,
resolving the unity:
$$
\frac1{(2\pi\hb)^n}\int_{\cX} \Pi^0_{y,\hb}\,dm(y)=1,
$$
and such that the rescaled functions
\begin{equation}
\pi^0_{y,\hb}(u)\od \Pi^0_{y,\hb}(\Exp_y(\sqrt{\hb}u)),\qquad 
u\in T_y\cX,
\tag{8.5}
\end{equation}
depend regularly on $\hb\to+0$.

We define the family of {\it dynamically deformed 
coherent states\/} $\Pi^t_{y,\hb}$ by applications of the
oscillator type evolution operators:
\begin{equation}
\wh{\Pi}^t_{y,\hb}\od
\underset{\leftarrow}{\Exp}
\bigg\{-\frac i\hb\int^t_0\wh{\bH}^t_y\,dt\bigg\}
\wh{\Pi}^0_{y,\hb}.
\tag{8.6}
\end{equation}
Using formula (8.4) and integrating over $y\in\cX$, 
we obtain a representation of the quantum evolution 
in the following form.

\begin{theorem}
Let the initial state {\rm(}Wigner function{\rm)} 
$\rho_\hb$ be resolved by the family of coherent states 
$\Pi^0_{y,\hb}$ as follows{\rm:}
$$
\rho_\hb(x)=\frac1{(2\pi\hb)^n}\int_{\cX}\Pi^0_{y,\hb}(x)
\rho^0_\hb(y)\,dm(y),
$$
where $\rho^0_\hb$ is a distribution over~$\cX$.
Then the Schr\"odinger type evolution of the quantum state 
$\hat{\rho}_\hb$ is given by the formula
\begin{equation}
\exp\bigg\{-\frac{it}{\hb}\wh{H}\bigg\}\hat{\rho}_\hb
=\frac1{(2\pi\hb)^n}\int_{\cX} \rho^0_\hb(y)
\exp\bigg\{-\frac{it}{\hb}H(y)\bigg\}
\wh{U}_{X^t(y)}\wh{\Pi}^t_{y,\hb}\,dm(y).
\tag{8.7}
\end{equation}
In this formula, the quantum Hamilton trajectory 
$\wh{U}_{X^t(y)}$ is determined by {\rm(7.5):}
\begin{equation}
\wh{U}_{X^t(y)}=\underset{\leftarrow}{\Exp}
\bigg\{-\frac{i}{2\hb}\int_{X^t(y)}\wh{\cH}^\hb\bigg\},
\tag{8.8}
\end{equation}
where $\wh{\cH}^\hb$ is the quantum Ether Hamiltonian
over~$\cX$, and the multiplicative integral 
in {\rm(8.8)} is taken along  
the Hamilton trajectory~{\rm(6.2)}.

The semiclassical asymptotics of the quantum trajectory 
{\rm(8.8)} is determined by Theorem~{\rm7.2:}
\begin{equation}
U_{X^t(y)}(x)
=\exp\bigg\{\frac{i}{\hb}\int_{\Sigma^t_y(x)}\omega\bigg\}
\varphi^t_y(x)+O(\hb),
\tag{8.9}
\end{equation}
where
\begin{equation}
\varphi^t_y(x)
=2^n\det\big(I-D(s_x\circ[X^t_y])(\tilde{x}^t_y)\big)^{-1/2}.
\tag{8.10}
\end{equation}
In {\rm(8.10)} we denote by $\tilde{x}^t_y$ 
the fixed point of the symplectic transformation 
$s_x\circ[X^t(y)]$, 
where $s_x$ are the reflection mappings related via 
{\rm(2.3), (2.4)} to the Ether Hamiltonian~$\cH$, 
and $[X^t(y)]$ is the symplectic trajectory~{\rm(4.1):}
\begin{equation}
[X^t_y]^*=\underset{\to}{\Exp}
\bigg(\frac12\int_{X^t(y)}\ad(\cH)\bigg).
\tag{8.11}
\end{equation}
The boundary of the membrane $\Sigma^t_y(x)$ in {\rm(8.9)}
consists of an arbitrary path~$c$ from~$\tilde{x}^t_y$ to~$y$,
the trajectory {\rm(6.2)} from~$y$ to~$X^t(y)$,
the path $[X^t(y)](c)$ 
{\rm(}with the opposite orientation{\rm)}, 
and the Ether geodesic connecting 
$s_x(\tilde{x}^t_y)$ with $\tilde{x}^t_y$
through the mid-point~$x$.
\end{theorem}

Formulas (8.6)--(8.10) generalize 
to the case of symplectic manifolds 
the well-known Gaussian type approximation program 
developed for Euclidean spaces 
(see, for instance, in \cite{Bab,Gui,MS,VPM,P,Kl,F1}).
We stress that in (8.7) the quantum translocation operator 
$\wh{U}_{X^t(y)}$ is separated from the deformed coherent 
states $\wh{\Pi}^t_{y,\hb}$.
Such a separation (or factorization) is a version of the
so-called ``interaction representation'' in quantum physics.
In our case the role of the leading free motion is played by
the intrinsic Ether stream.

The Gaussian type approximation is a particular case of the
general scheme described above.

Namely, let $\{\lambda^0_y\}$ be a positive Lagrangian
distribution in the complexified tangent bundle ${}^{\bC}T\cX$
(say, given by an almost complex structure on~$\cX$).
Denote by $\{\lambda^t_y\}$ a distribution obtained 
from $\{\lambda^0_y\}$ by rotating each Lagrangian
plane~$\lambda^0_y$ by means of the linearized Hamiltonian
system (6.9).
Let the functions $\pi^0_{y,\hb}$ at $\hb=0$ be just Gaussian
exponents assigned to $\lambda^0_y$ 
as in \cite{Gui,MS}.
Then the dynamically deformed coherent states 
$\Pi^t_{y,\hb}$ (8.6)
can also be presented in the form (8.5): 
\begin{equation}
\Pi^t_{y,\hb}(\Exp_y(\sqrt{\hb}u))\od \pi^t_{y,\hb}(u),
\tag{8.12}
\end{equation}
where $\pi^t_{y,\hb}$ are regular in $\hb\to+0$ 
and $\pi^t_{y,0}$ are just Gaussian exponents assigned 
to~$\lambda^t_y$.

The function $\pi^t_{y,0}$ is evaluated explicitly by solving
the Cauchy problem for the harmonic oscillator type equation
over the tangent space $T_y\cX$ as follows:
\begin{equation}
\pi^t_{y,0}(u)=\underset{\leftarrow}{\Exp}
\bigg\{-i\int^t_0 Q^\tau_y
(u-i\Psi(y)\pa_u)\,d\tau\Big\}\pi^0_{y,0}(u).
\tag{8.13}
\end{equation}
Here the oscillator type Hamiltonian is given by 
$$
Q^t_y(u)\od \frac12\langle V^{t*}_y \nabla^2 H(X^t(y))
V^t_y u,u\rangle,\qquad u\in T_y\cX.
$$ 

For example, if the connection $\Gamma$ on $\cX$ is chosen 
so that the function~$H$ is covariantly quadratic (see (6.11)),
then (8.13) reads
\begin{equation}
\pi^t_{y,0}(u)=\exp\bigg\{-\frac{it}{2}
\langle\nabla^2 H(y)(u-i\Psi(y)\pa_u),
u-i\Psi(y)\pa_u\rangle\bigg\}\pi^0_{y,0}(u).
\tag{8.14}
\end{equation}
Since $\pi^0_{y,0}$ is a Gaussian exponent, the application of
the oscillator-type evolution operators in (8.13), (8.14) is
made easily and explicitly by the well-known formulas. 

\setcounter{section}{1}
\setcounter{theorem}{0}
\renewcommand{\thesection}{\Alph{section}}

\section*{Appendix.\\ Extension of Cartan's moving frame method}

In this Appendix we represent a nonsymplectic version 
of the dynamic geometry from Sects.~2--6
making this theory to be applicable to arbitrary manifolds. 

Let $\cM$ be a manifold with an affine connection $\Gamma$. 
We denote by $\nabla$ the covariant derivative with respect to
$\Gamma$, i.e.,  
$$
(\nabla_{\!j}\,\, u)^k=\pa_j u^k+\Gamma^k_{sj} u^s
$$
for any vector field $u$ on $\cM$. 

Consider a differential one-form $A$ on $\cM$ with values in the
space of vector fields on $\cM$.  
This form can also be presented as a family of linear mappings 
$A_x(z):T_x\cM\to T_z\cM$ (where $x,z\in\cM$) so that 
the form~$A$ at the point~$x$ is given by 
$$
A_x=A_x(\cdot)_j\, dx^j,\qquad A_x(z)_j\in T_z\cM.
$$
We denote by $a=A|_{\diag}$ the diagonal family of mappings
$T_x\cM\to T_x\cM$ (where~$x$ is running over~$\cM$). 

Assume that 
\begin{equation}
\pa A+\frac12[A\wedgeco A]=0,
\tag{A.1}
\end{equation}
where $\pa$ is the differential (acting in the space of forms on
$\cM$) and $[\cdot,\cdot]$ is the commutator (acting in the
space of vector fields on $\cM$).
This is an analog of Eq.~(2.1).

The field $A$ generates a connection $\pa+A$ 
in the trivial bundle over~$\cM$ with the fiber $C^\infty(\cM)$.
Condition (A.1) means that this connection is curvature free.
Therefore, we obtain a big family of {\it internal translations\/}
$g_{x,y}:\,\cM\to\cM$ labeled by pairs of points $x,y,\in\cM$.
By definition, 
$$
g^*_{x,y}=\underset{\to}{\Exp} \bigg(\int^x_y A\bigg).
$$
The diffeomorphisms $g_{x,y}$ are analogs of 
Ether transformations (see Sect.~2).
Their definition does not depend on a path connecting the
points~$y$ and~$x$ because of the zero curvature
condition~(A.1). 

Note that solutions~$A$ of Eq.~(A.1)
are essentially different from solutions
of a zero curvature equation used in~\cite{EmW}
following the deformation quantization approach~\cite{F}.
The geometric objects used in \cite{EmW,F} are
fiberwise (vertical) fields on $T\cM$.
In particular, 
the parallel translation from~$y$ to~$x$
generated by the zero curvature connection 
in the sense of~\cite{EmW} 
is a mapping $T_{y}\cM\to T_{x}\cM$
but not a diffeomorphism of $\cM$ as in our approach.
This is exactly the difference between kinematics and dynamics
as it was explained in the Introduction. 

Besides the zero curvature equation (A.1), 
we fix the boundary condition
\begin{equation}
\nabla' A \Big|_{\diag}=0,
\tag{A.2}
\end{equation}
where $\nabla'$ denotes the covariant derivative of a
vector-valued form with respect to the adjoint
connection~$\Gamma'$, i.e., 
$(\nabla'_k A)^s_j\od D_k A^s_j+\Gamma^s_{kl}A^l_j$.

Combining (A.1) and (A.2), 
we obtain the following equation on the diagonal:
\begin{equation}
\delta a+\frac12( a\wedgeco a)=0.
\tag{A.3}
\end{equation}
Here $\delta$ denotes the adjoint covariant differential
$\delta=d+\Gamma'\wedge$.
The brackets $(\cdot,\cdot)$ are generated by the {\it  
torsion tensor\/} $T$ of the connection~$\Gamma$:
$$
(u,v)^k \od u^s T^k_{sl}v^l.
$$

Equation (A.3) coincides with the first structural equation 
in the moving frame method due to E.~Cartan.
The solvability of (A.3) is guaranteed by the 1-st Bianchi's
identity. 
This structural equation combines the covariant strength 
$\delta a$ of the diagonal vector field~$a$ with the torsion of
the connection~$\Gamma$
(see \cite{OR} about Cartan's ideas concerning the torsion). 

Thus we see that the zero curvature equation (A.1) is a natural
extension of the structural equation from~$\cM$ to 
$\cM\times\cM$, and the field~$A$ is an extension 
of the {\it Cartan field}~$a$. 

By choosing a field~$a$,
we set in addition to (A.2), the boundary condition 
\begin{equation}
A\Big|_{\diag}=a.
\tag{A.4}
\end{equation}

Any solution of (A.1) satisfying the boundary conditions (A.2)
and (A.4) will be called an {\it internal vector\/}
field on~$\cM$. 

Let us consider ``trajectories'' of the internal vector
field~$A$, that is, solutions of the equation
\begin{equation}
\pa s=A(s).
\tag{A.5}
\end{equation}
Here $\pa=\pa/\pa x$ plays the role of time derivative, and
the ``trajectory'' $s=s_x(z)$ of (A.5) is uniquely determined by
the initial data
$$
s_x(z)\Big|_{x=z}=z.
$$
This is the analog of the Hamiltonian dynamics (2.4).
Condition (A.1) guarantees the solvability of (A.5).

Assuming that the field~$A$ is complete, 
we obtain a family $\{s_x\mid x\in\cM\}$ of diffeomorphisms
of the manifold~$\cM$. For each $x\in\cM$ the mapping $s_x$ 
has the fixed point $z=x$, i.e., 
\begin{equation}
s_x(x)=x.
\tag{A.6}
\end{equation}
Obviously, in the domain where $\det a(x)\ne0$, 
this fixed point is isolated.

We shall say that the manifold $\cM$ is endowed with 
an {\it inversive structure\/} if there is a family of
diffeomorphisms $\{s_x\mid x\in\cM\}$ possessing isolated fixed
points (A.6). 
In general, the inversions $s_x$ are not involutions 
$(s^2_x\ne\id)$ and so the inversive structure is not 
a reflective structure.

\begin{proposition}
{\rm(i)} Any inversive structure generates 
an internal vector field 
by the formula
$$
A_x(z)=(\pa_x s_x)(s^{-1}_x(z)).
$$
The corresponding affine connection and the Cartan field
are given by 
\begin{equation}
\Gamma^l_{jk}(x)=-\frac{\pa^2 s^l_x(z)}{\pa z^m \pa x^r}
\bigg[\frac{\pa s_x(z)}{\pa z}\bigg]^{-1\,m}_{\,\,k}
\bigg[\frac{\pa s_x(z)}{\pa x}\bigg]^{-1\,r}_{\,\,j}
\bigg|_{z=x},\qquad
a(x)^l_j=\frac{\pa s^l_x(z)}{\pa x^j}
\bigg|_{z=x}.
\tag{A.7}
\end{equation}
The Cartan field~$a$ has no eigenvalues~$0$ or~$1$ in the spectrum.

{\rm(ii)}
Let $A=A^+$ be an interval vector field on $\cM$. 
Let $s^+$ be the family of inversions generated by~$A^+$ via
{\rm(A.5)}. 
Then there exists another solution of the zero curvature
equation {\rm(A.1):}
\begin{equation}
A^- \od -(Ds^+)^{-1}\cdot A^+(s^+),
\tag{A.8}
\end{equation}
or in more detail, 
$$
A^-_x(z)=-\bigg[\frac{\pa s^+_x(z)}{\pa z}\bigg]^{-1}
\cdot A^+_x(s^+_x(z)).
$$
The inversions $s^-$ generated by $A^-$ via {\rm(A.5)} 
are just $s^-_x=(s^+_x)^{-1}$.
The boundary conditions are
$$
A^-\bigg|_{\diag}=a^-,\qquad (\nabla^-)' A^-\bigg|_{\diag}=0.
$$
Here $a^-=a^+/(a^+-I)$, 
and the covariant derivative $\nabla^-$ is taken with respect 
to the connection $\Gamma^-$ defined by {\rm(A.7)} 
via the inversions~$s^-$.

The diagonal field $a^-$ is a Cartan field, that is, 
the solution of {\rm(A.3)} corresponding to the
connection~$\Gamma^-$.  

{\rm(iii)}
The internal translations $g_{x,y}$ are related to inversions 
as follows {\rm:}
$$
g_{x,y}=s^+_x\circ s^-_y.
$$

{\rm(iv)}
Let $E^\pm=\Exp^\pm_x(vt)$ denote the solutions of the equations
$$
\frac{d}{dt}E^\pm=\frac12\langle A^\pm_x(E^\pm),v\rangle,\qquad
E^\pm\bigg|_{t=0}=x, 
$$
where $v\in T_x\cX$. Then
\begin{equation}
s^+_x(\Exp^-_x(v))=\Exp^+_x(-v).
\tag{A.9}
\end{equation}
\end{proposition}

The curve $\sigma^+\cup \sigma^-\subset \cM$ 
composed of two pieces 
$\sigma^+=\{\Exp^+_x(vt)\mid t\geq0\}$ and 
$\sigma^-=\{\Exp^-_x(vt)\mid t\leq0\}$ 
we call an {\it internal geodesic through the center}~$x$. 
It follows from (A.9) that such internal geodesics 
are {\it inversive curves\/} with respect to the
inversion~$s^\pm$. 

The pair of vector fields~$A^\pm$ described in
Proposition~A.1,\,(ii)
we call an {\it internal pair}.

In the case of a symplectic manifold~$\cM$ 
(with a symplectic form~$\omega$),
the internal vector fields~$A^\pm$ are Hamiltonian:
\begin{equation}
A^\pm_x(z)=D\cH^\pm_x(z)\Psi(z),
\qquad\text{where}\quad
\Psi=\omega^{-1}\quad\text{is a Poisson tensor}.
\tag{A.10}
\end{equation}
The inversions $s^\pm_s$ defined by the Hamiltonian
systems~(A.5) are symplectomorphisms of~$\cM$.

\begin{proposition}
In the symplectic case 
the connection $\Gamma=\Gamma^\pm$ defined by {\rm(A.7)}
via inversions $s=s^\pm$
is a symplectic connection{\rm:}
$$
\nabla^\pm \omega=0.
$$
The torsion $T=T^\pm$ of this connection obeys 
the cyclicity condition 
\begin{equation}
\underset{j,k,l}{\mathfrak S}\omega_{js} T^s_{kl}=0
\qquad
(\text{summation over cyclic permutations}).
\tag{A.11}
\end{equation}
The internal Hamiltonian $\cH=\cH^\pm$ {\rm(A.10)} 
is reconstructed from inversions $s=s^\pm$ by the formula
\begin{equation}
\cH_x(z)=\int^z_x\langle \pa s_x(s^{-1}_x(z)),\omega(z)\,dz\rangle.
\tag{A.12}
\end{equation}
Here we assume that the zero boundary conditions hold 
on the diagonal{\rm:} 
\begin{equation}
\cH^\pm_x(z)\bigg|_{x=z}=0.
\tag{A.13}
\end{equation}
Two internal Hamiltonians $\cH^+$ and $\cH^-$ 
are related to each other by the inversion mapping{\rm:}
\begin{equation}
\cH^+_x(s^+_x(z))=-\cH^-_x(z).
\tag{A.14}
\end{equation}
The derivatives of $\cH=\cH^\pm$ satisfy the following boundary
condition on the diagonal{\rm:}
\begin{equation}
\nabla_l\nabla_m\cH_k\bigg|_{\diag}
=\omega_{ms}T^s_{lr}\Psi^{rj}\nabla_j\cH_k\bigg|_{\diag}.
\tag{A.15}
\end{equation}
Here the covariant derivatives $\nabla=\nabla^\pm$ 
and the torsion tensor $T=T^\pm$ 
are assigned to the connection $\Gamma=\Gamma^\pm$.
\end{proposition}

So we see that, in the symplectic case,
formula~(A.12) replaces formula~(2.7), 
relation~(A.14) replaces the skew-symmetry relation~(2.3), 
the boundary conditions~(A.13), (A.14) replace conditions~(2.2).

Note that the general situation 
where we deal with torsion and the inversive structure 
is, in fact, very common and presented in many examples.
Nevertheless, let us now consider the torsion free case 
and the reflective structures on a manifold~$\cM$.

\begin{proposition}
{\rm(i)}
The inversions $s_x$ defined by {\rm(A.5)} 
are involutions 
iff the following skew-symmetry condition holds{\rm:}
\begin{equation}
A(s)=-Ds\cdot A,
\tag{A.16}
\end{equation}
or in detailed notation{\rm:}
$$
A_x(s_x(z))=-\frac{\pa s_x(z)}{\pa z} A_x(z).
$$
In this case, $A=A^+=A^-$,
the boundary condition holds{\rm:}
\begin{equation}
a=A\Big|_{\diag}=2,
\tag{A.17}
\end{equation}
the connection $\Gamma$ is given by {\rm(2.6)},
and $\Gamma$ is torsion free{\rm:} $T=0$.

{\rm(ii)}
For any torsion free affine connection on~$\cM$, 
the zero curvature equation {\rm(A.1)} with boundary
conditions {\rm(A.2), (A.17)} and with 
the skew-symmetry condition {\rm(A.16)} 
has a solution $A=A_x(z)$ at least near the
diagonal $\diag=\{x=z\}$ in $\cM\times\cM$,
and in this case the {\rm(}semiglobal{\rm)}
reflective structure is given on~$\cM$.  
\end{proposition}

An internal vector field $A$ on $\cM$ satisfying conditions 
(A.16), (A.17) will be called a {\it fundamental vector field}.

Using a fundamental vector field on $\cM$ and
multiplying it by the factor~$\frac12$, we obtain analogs of the
results in Sects.~3--6. Indeed,
it is clear that Hamiltonian systems (2.5), (3.1), (3.2) are
easily extended to the case of general 
fundamental vector field. 
In particular, 
the notion of path-diffeomorphism is well defined
and the analog of formula (4.1) holds:
$$
[\sigma]^*=\underset{\to}{\Exp}\bigg(\frac12\int_\sigma A\bigg).
$$

This is the parallel translation along the path~$\sigma$ 
with respect to the connection 
$$
{\boldsymbol\nabla}^0=\pa+\frac12 A
$$
acting in the trivial bundle over $\cM$ 
with the fiber $C^\infty(\cM)$.

Theorems~2.2, 3.1, 4.2 hold in this general case as well.
Instead of the Ether curvature {\rm(5.3)}, 
there appear the {\it dynamic curvature\/} vector fields
$$
B_{jk}=\frac14[A_{0k},A_{0j}]
=[{\boldsymbol\nabla}^0_j,{\boldsymbol\nabla}^0_k]
$$
generating the Lie algebra of the dynamic holonomy group
$\cL_0$ of the manifold~$\cM$. 
The analogs of formulas in Theorem~5.1 are the following:
$$
B^s_{jk}(0)=0,\qquad 
(\nabla_{\!m}\,\, B_{jk})^s (0)=2 R^s_{mjk}(0).
$$
Here the covariant derivative~$\nabla$ 
corresponding to the connection~$\Gamma$.

The {\it dynamic curvature form\/} is defined as 
$$
B=\frac12 B_{jk}dx^k\wedge dx^j =\frac18[A\wedgeco A].
$$
This is the $2$-form on $\cM$ with values in vector fields on
$\cM$. The analog of formula (5.1) is 
$$
[\pa\Sigma]^{-1*}=\exp\bigg(\int_{\Sigma}B+o(\Sigma)\bigg).
$$

In the presence of a general fundamental field 
the translocation operation is described 
as follows. 
We start from a dynamical system 
\begin{equation}
\frac{d}{dt}X=u(X),\qquad X\bigg|_{t=0}=x,
\tag{A.18}
\end{equation}
where $u$ is a vector field on~$\cM$. 
The translation along 
trajectories of (A.18) is denoted by $X=X^t(x)$. 

The segment of the trajectory (6.2) generates the
path-diffeomorphism $[X^t(y)]:\,\cM\to\cM$.
Here the point $y\in\cM$ is arbitrarily fixed.

Then we define a time-dependent vector field on~$\cM$:
\begin{equation}
v^t_y\od[X^t(y)]^{-1}_*\,\Big(u-\frac12 \dot X^t(y) A_{X^t(y)}\Big).
\tag{A.19}
\end{equation}
Here $A$ is the fundamental field on~$\cM$, and the
subscript~$*$ means the standard operation on vector fields:
$$
(\gamma^{-1}_* u)(z)\od [d\gamma(z)]^{-1}\,u(\gamma(z)),
\qquad \gamma:\,\cM\to\cM.
$$

The field (A.19) determines the new dynamical system on~$\cM$:
\begin{equation}
\frac{d}{dt}Z=v^t_y(Z),\qquad Z\bigg|_{t=0}=x.
\tag{A.20}
\end{equation}
The solution of this system we denote by $Z=Z^t_y(x)$.

The transformation from (A.18) to (A.20) is the {\it
translocation\/} to~$y$. 
There is the following analog of Theorem~6.1.

\begin{proposition}
{\rm(i)} Solutions of the original system {\rm(A.18)} 
and the translocated system {\rm(A.20)} 
are related to each other by the formula 
$$
X^t=[X^t(y)]\circ Z^t_y.
$$
Here $[X^t(y)]$ is the path-diffeomorphism corresponding 
to the trajectory of {\rm(A.18)} starting at~$y$.

{\rm(ii)} The point $y$ is an equilibrium point 
of the translocated-to-$y$ system {\rm(A.20)}.
One has
$$
v^t_y(y)=0,\qquad (\nabla v^t_y)(y)=M^t_y,
$$
where
\begin{equation}
M^t_y=(V^t_y)^{-1}\cdot \nabla u(X^t(y))\cdot V^t_y.
\tag{A.21}
\end{equation}
Here the covariant derivative $\nabla$ corresponds 
to the connection~$\Gamma$. 

{\rm(iii)} The differential of the solution to system {\rm(A.18)}
at the point~$y$ is given by formulas {\rm(6.8), (6.9)},
where the linear maps $M^t_y:\, T_y\cM\to T_y \cM$ are defined
by~{\rm(A.21)}. 
\end{proposition}

\begin{corollary}
If the vector field~$u$ and the
connection~$\Gamma$ satisfy the consistency condition
\begin{equation}
\nabla_u(\nabla u)=0
\tag{A.22}
\end{equation}
along the trajectory $X^t(y)$ {\rm(A.18)} starting at~$y$, 
then
\begin{equation}
d X^t(y)=V^t_y \cdot \exp(t\nabla u(y)).
\tag{A.23}
\end{equation}
For the periodic trajectory, 
formula {\rm(A.23)}
means the factorization of the monodromy matrix 
to the geometric monodromy
{\rm(}the holonomy of~$\Gamma${\rm)}
and the dynamic monodromy
{\rm(}with generator $\nabla u(y)${\rm)}.
\end{corollary}

If condition (A.22) holds at any point, 
then we say that the vector field ~$u$ is 
{\it auto-linear\/} with respect to the connection~$\Gamma$,
or that the connection~$\Gamma$ {\it covariantly linearizes\/}
the field~$u$. 

In respect of the factorization formula (A.23), 
the following natural questions arise:
How a connection which covariantly linearizes a given
vector field can be found?
How a connection which covariantly linearizes a given
vector field along a given trajectory can be found?

\bibliographystyle{amsalpha}

\begin{thebibliography}{99}

\bibitem{K}
M.~Karasev, 
{\it Quantization and intrinsic dynamics}.
In: M.~Karasev (ed.), 
{\it Asymptotic Methods for Wave and Quantum Problems}, 
AMS Trans., Ser.~2, Vol.~208, pp.~1--32,
AMS, Providence, 2003.
Preprint version in arXiv: math.QA/0207047.

\bibitem{K1}
M.~Karasev, 
{\it Intrinsic dynamics of symplectic manifolds: 
membrane representation and phase product},
Preprint version in arXiv: math.QA/0308118.

\bibitem{KM}
M.~Karasev and V.~Maslov,
{\it Nonlinear Poisson brackets. Geometry and Quantization}, 
Nauka, Moscow, 1991; 
English transl., Amer. Math. Soc., Providence, RI, 1993.

\bibitem{Bab}
V.~M.~Babich, 
{\it The Multidimensional WKB Method and the Ray method.
Its Analogues and Generalizations},
Encyclopedia of Math. Sci., Vol.~34 (Partial Diff. Eq.,~V)
Springer-Verlag, 1993.

\bibitem{MRR}
J.~Marsden, T.~Ratiu, and G.~Raugel,
{\it Symplectic connections and the linearization of Hamiltonian
systems}, Proc. Roy. Soc. Edinburgh Ser.~A, {\bf 117} (1991),
329--380. 

\bibitem{K2}
V.~P.~Maslov,
{\it Operator Methods},
Nauka, Moscow, 1973 
(in Russian).

\bibitem{Gui}
V.~Guillemin,
{\it Symplectic Spinors and Partial Differential Equations},
Colloques Intern. C. N. R. S., 1975, v.~237.

\bibitem{MS}
A.~Melin and J.~Sj\"ostrand,
{\it Fourier integral operators with complex-valued phase
fucntions}, 
Lect. Notes Math., {\bf 459} (1975), 121--223.

\bibitem{VPM}
V.~P.~Maslov,
{\it The Complex WKB-Method in Nonlinear Equations},
Nauka, Moscow, 1977 (in Russian);
English transl., 
{\it The Complex WKB-Method in Nonlinear Equations},
I. Linear Theory, Birkh\"auser, Basel--Boston--Berlin, 
1994. 

\bibitem{P}
M.~Popov, 
{\it A new method for calculating wave fields 
in high frequency approximation},
Zapiski Nauch. Sem. LOMI, {\bf 104} (1981), 195--216
(in Russian).

\bibitem{Kl}
J.~Klauder,
{\it Global uniform semiclassical approximation to wave
equation}, 
Phys. Rev. Lett., {\bf 56}(9) (1986), 897--899.

\bibitem{F1}
B.~Fedosov, 
{\it A trace formula for Schr\"odinger operator},
Russ. J. Math. Phys., {\bf 1}(4) (1993), 447--463.

\bibitem{OR}
L.~O'Raifeartaigh, 
{\it The Dawning of Gauge Theory}, Princetion Univ. Press, 
New Jersey, 1997. 

\bibitem{EmW}
C.~Emmrich and A.~Weinstein,
\textit{The differential geometry of Fedosov's quantization}, 
in: \textit{Lie Theory and Geometry. In Honor of B.~Kostant},
Progr. Math. \textbf{123} (1994), Birkh\"auser, New York, 217--240.

\bibitem{F}
B.~Fedosov,
\textit{A simple geometrical construction of deformation
quantization},
J. Diff. Geom. \textbf{40} (1994), 213--238.


\end{thebibliography}

\end{document}